\documentclass[11pt,twoside,reqno]{amsart}
\usepackage{epsfig}

\calclayout

\newtheorem{theorem}{\bf Theorem}

\newtheorem{proposition}[theorem]{\bf Proposition}

\newtheorem{defn}{\bf Definition}

\newtheorem{con}{\bf Conjecture}
\usepackage{float}
\usepackage{calc,tikz,amsmath}
\usetikzlibrary{arrows,decorations.markings,positioning,shapes.geometric,}
\tikzstyle{vertex}=[circle, draw, inner sep=1pt, minimum size=4pt]
\newcommand{\vertex}{\node[vertex]}

\tikzstyle{ann} = [fill=white,font=\footnotesize,inner sep=1pt]
\tikzstyle{arrow} = [thick,<-->,>=stealth]

\begin{document}

%\leftline{ \scriptsize \it  Journal of Prime Research in Mathematics
%Vol. {\bf 17}(2) (2021), 1-7}

\vspace{1.3cm}

\title
{Parity Labeling in Signed Graphs}

\author{Mukti Acharya$^1$, Joseph Varghese Kureethara$^2$}
\thanks{ {\enskip
  \enskip $^{1}$CHRIST (Deemed to be University), Bangalore 560029, Karnataka, India. Email: mukti1948@gmail.com\\
$^{2}$CHRIST (Deemed to be University), Bangalore 560029, Karnataka, India. Email: frjoseph@christuniversity.in}}
\begin{abstract}
Let $S=(G, \sigma)$ be a signed graph where $G=(V, E)$ is a graph called the underlying graph of $S$ and $\sigma:E(G) \rightarrow \{+,~-\}$. Let $f:V(G) \rightarrow \{1,2,\dots,|V(G)|\}$ such that $\sigma(uv)=+$ if and only if $f(u)$ and $f(v)$ are of same parity and $\sigma(uv)=-$ if and only if $f(u)$ and $f(v)$ are of opposite parity. Under $f$ we get a signed graph $G_f$ denoted as $S$, which is a parity signed graph. In this paper, we initiate the study of parity labeling in signed graphs and we define and find `rna' number denoted as $\sigma^-(S)$ for some classes of signed graphs. We also characterize some signed graphs which are parity signed graphs. Some directions for further research are also suggested. \vskip 0.4 true cm
 \noindent
 {\it AMS SUBJECT} : Primary 05C22, 05C75.\\
\end{abstract}
\maketitle

%\TagsOnRight

\pagestyle{myheadings} \markboth{\centerline {\scriptsize Acharya,
Kureethara}}
         {\centerline {\scriptsize  Parity labeling in Signed Graphs }}

%%%%%%%%%% the following introduction form is an option %%%%%%%%%%%%

\bigskip
\bigskip
\medskip

\section{Introduction}

The concept of signed graph is very popular in Graph Theory. Here we introduce a type of signed graphs called as \textit{parity signed graphs}. This is based on the assignment of the $+$ and $-$ sign of its edges of a graph based on the positive integer labels to the vertices of a graph. For terminologies of graphs we refer to \cite{H69,W} and for signed graphs we refer to \cite{Z82}. All graphs/signed graphs considered here are simple and connected unless mentioned otherwise. There is a tremendous growth on the study of signed graphs as there are numerous applications for it in both industrial and theoretical realms. For a detailed conceptual framework in signed graphs, we refer the reader to \cite{Z05}.

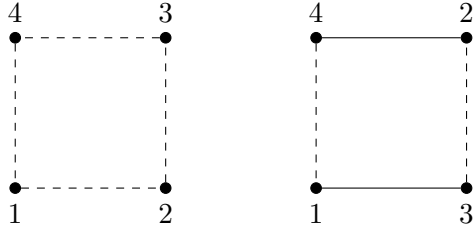
\begin{figure}
	\centering
	\begin{tikzpicture}%[scale=0.5]
	%cycle C_4
	\vertex (1) at (0,0)[label=below:$1$] [fill]{};
	\vertex (2) at (2,0)[label=below:$2 $] [fill]{};
	\vertex (3) at (2,2) [label=above:$3$][fill]{};
	\vertex (4) at (0,2)[label=$4$] [fill]{};
	%second cycle
	\vertex (5) at (4,0)[label=below:$1$] [fill]{};
	\vertex (6) at (6,0)[label=below:$3 $] [fill]{};
	\vertex (7) at (6,2) [label=above:$2$][fill]{};
	\vertex (8) at (4,2)[label=$4$] [fill]{};
	\path 
	%cycle C_3
	(1) edge [dashed](2)
	(2) edge [dashed](3)  
	(3) edge [dashed](4)
	(4) edge [dashed](1)
	%second cycle
	(5) edge (6)
	(6) edge [dashed](7)  
	(7) edge (8)
	(8) edge [dashed](5)  
	;
	\end{tikzpicture}
	\caption{Two distinct parity labelings of cycle $C_4 $.}\label{fig:pscycle}
	
\end{figure}

By an $(n,m)$-graph, $G=(V,E)$, we mean a graph $G$ such that $n=|V(G)|$ and $m=|E(G)|$. A signed graph $S=(G, \sigma)$ is such that $G=(V,E)$ and $\sigma:E(G)\rightarrow \{+,-\}$ and edges which receive + ($-$) signs are called positive (negative) edges of $S$, respectively. The graph $G$ is called the underlying graph of $S$. By $E^+(S)(E^-(S))$ we denote the set of positive (negative) edges of $S$ and the edge set $E(S)=E^+(S)\cup E^-(S)$. A signed graph is said to be all-positive if $E^-(S)=\emptyset$ and all-negative if $E^+(S)=\emptyset$. While drawing the signed graph, positive edges are drawn as solid line segments and negative edges as dashed line segments as depicted in the FIGURE \ref{fig:pscycle}. A signed graph is said to be homogeneous if it is either all-positive or all-negative and heterogeneous otherwise. Here, by a positive (negative) homogeneous signed graph we mean a signed graph which is all-positive (all-negative) \cite{AJK}. In this paper, we initiate the study of parity labeling in signed graphs which are known as parity signed graphs.
%pic-1
\begin{defn}
	A signed graph $S$=(G, $\sigma$) is a \textbf{parity signed graph}, if there exists a bijection $f:V(G) \rightarrow \{1, 2, \dots, n\}$ such that for an edge $uv$ in $G$, $\sigma(uv)=+$ if $f(u)$ and $f(v)$ are of same parity and $\sigma(uv)=-$ if $f(u)$ and $f(v)$ are of opposite parity.
\end{defn}

There are many different types of signed structures for any given graph $G$. It depends on the bijection $f$. In that context, the signed graph $S$ is better expressed as $G_f$. However, we use the notation $S$ unless there is particular requirement of distinctly identifying two parity signed graphs generated from the same graph $G$.

A parity signed graph can be homogeneous or heterogeneous. For example, given a path, P${_n}$ with $n$ vertices, we can see that every edge can take the `$-$' sign, that is, we have all-negative paths P$_n$ by assigning the vertices the labels 1 to $n$ such that every pair of adjacent vertices gets consecutive integers from 1 to $n$. In this labeling, adjacent vertices receive opposite parity labeling. Thus negative homogeneous P$_n$ is a parity signed graph.

It is not difficult to see that every graph has a parity signed labeling. However, every signed graph need not be a parity signed graph. This is because, the labeling need not be parity signed labeling. For example, let us consider the cycle C$_3$. As the vertex labels are 1, 2 and 3, no parity signed labeling exists for homogeneous C$_3$. This is because of the fact that every integer, odd or even, among 1, 2 and 3 is adjacent to an odd integer and an even integer. We could extend this observation to all complete graphs of order at least 3.

The above discussion leads us to search for those signed graphs that are parity signed graphs.  Before going deeper into finding such characterizations, we inspect the signed structure of some standard graphs. We had already seen that an all-negative path is a parity signed graph.

Is path the only negative homogeneous parity signed graph? The answer is no. There are cycles that are negative homogeneous parity signed graph.  Throughout the text, by an odd (even) vertex, we mean a vertex that receives odd (even) integer as its label.

\begin{theorem}\label{negc}
	A negative homogeneous signed cycle of order $n\ge4$ is a parity signed graph if and only if $n$ is even. 
\end{theorem} 
\begin{proof}
	Let a negative homogeneous signed cycle of order $n\ge4$ be a parity signed graph. It means that end vertices of every edge receive labels of different parity. However, a cycle has equal number of vertices and edges. As every edge has its ends labeled by integers of opposite parity, there are equal number of odd and even integers. This implies that the total number of edges is an even number. We conclude that the total numbers of vertices is an even number.
	
	Assume that $n$ is even. Then there are $\frac{n}{2}$ odd and even integers. Arranging them in a consecutive way on the vertices of a signed graph, we get the parity signed labeling.
\end{proof}

An immediate effect of this result is that if a graph has an odd cycle, then it is not a negative homogeneous parity signed graph.

It is easy to see that an all-positive path is not a parity signed  graph. As the labeling of the vertices begins with 1, we are forced to use the label 2 for some vertex. Hence, 2 must be adjacent to one of the odd numbers $\le n$, where $n$ is the number of vertices of the path. Hence, there is no parity signed path having all positive edges.

Is there any positive homogeneous parity signed graph? We answer this in negative in the following theorem.
\begin{theorem}\label{nont}
	Every connected non-trivial parity signed graph will have at least one negative edge.
\end{theorem}
\begin{proof}
	If $S$ is a connected non-trivial signed graph, then it has at least two vertices. By the same argument as in the case  of path we have seen that a vertex labeled with an even integer must be adjacent to a vertex labeled with an odd integer. Hence, there does not exist a connected non-trivial postive homogeneous parity signed graph.
\end{proof}

Since, every graph can have a labeling that gives a parity signed labeled graph and it is not immediate that a subsignedgraph of a parity signed graph is a parity signed graph. The labeling of a subgraph will be preserved in the subgraph, if the subgraph is formed by deletion of edges. However, for a connected graph, if a subgraph is formed by the deletion of a vertex labeled with an integer between 1 and $n$, the subgraph need not retain its signed structure. This is because of the possible relabeling of the vertices. These observations lead us to the following result. 

\begin{theorem}
	A connected subsignedgraph of a homogeneous parity signed graph need not be a parity signed graph.
\end{theorem}
\section{Optimization of the signd structure}
The search now is to find the smallest number of negative edges we can have in a parity signed graph. For any given graph $G$ with $n$ vertices, each bijection $f:V(G) \rightarrow \{1, 2, \dots, n\}$ gives a parity signed graph, not all may be distinct. In fact, $f$ is a permutation of the vertices of $G$. Hence, there are $n!$ such bijections. Taking into consideration the different parity signed structures of a graph, we define two parameters viz., the `rna' number and the `adhika' number of a parity signed graph.

The \textit{Sanskrit} word for `$-$' is `rna' means debt. It is pronounced as rina, where \textbf{ri} is as \textbf{ri} in \textbf{ri}bbon and \textbf{na} is \textbf{na} as in coro\textbf{na}. We define the `rna' number of signed graph as follows.
\begin{defn}
	The `rna' number of a parity signed graph $S$ is the smallest number of negative edges among all the parity signed labeling of its underlying graph $G$.
\end{defn}
The `rna' number is denoted by the symbol $\sigma^-(S)$, where $\sigma$ stands for word \textit{signum} means sign in Latin. We also define  the `adhika' number of a parity signed graph and it is denoted as $\sigma^+(S)$.
\begin{defn}
	The `adhika' number of a parity signed graph $S$ is the largest number of postive edges among all the parity signed labeling of its underlying graph $G$.
\end{defn}

The following inequalities are immediate for a parity signed graph $S$.
\begin{center}
	$\sigma^-(S) \le |E^-(S)|$ and $|E^+(S)| \le \sigma^+(S)$
\end{center}

Now, we find the `rna' number of some signed graphs, viz., paths, cycles, stars and complete graphs.

\begin{proposition}
	For any path P$_n$ with $n$ vertices, $\sigma^-(P_n)$=1.
\end{proposition}
\begin{proof}
	Let the vertices of the path be $v_i, i\in \{1, 2, \dots, n\}$ such that $v_iv_{i+1}\in E(P_n)$ for every $1\le i \le n-1$.
	
	Let $f:V(P_n)\rightarrow \{1, 2, \dots, n\}$ be the vertex labeling function such that 
	
	$f(v_i) = \begin{cases}
	2i-1 &\mbox{if } 1\le i \le \lceil{\frac{n}{2}}\rceil\\
	2i-(n+1) &\mbox{if } \lceil{\frac{n}{2}}\rceil+1 \le i \le n~and~n~is~ odd\\
	2i-n &\mbox{if } \lceil{\frac{n}{2}}\rceil+1 \le i \le n~and~n~is~even 
	\end{cases}$
	
	This labeling gives all the odd integers at first and then the even integers, consecutively to be vertices of the path. It is easy to see that the end vertices of exactly one edge get integers of different parity. Hence, $\sigma^-(P_n)$=1.
\end{proof}

\begin{proposition}\label{cyc2}
	For any cycle C$_n$ with $n$ vertices, $\sigma^-(C_n)$=2.
\end{proposition}
\begin{proof}
	We first label the vertices consecutively with odd integers from the set $A = \{1, 2, \dots, n\}$ and then with even integers. Observe that when the labeling switches from odd to even and then from the last even integer to the first odd integer (i.e.,  $n$ or $n-1$ to 1) two edges receive `$-$' sign.
	
	Let the vertices of the cycle be $v_i, i\in {1, 2, \dots, n}$.
	
	Let $f:V(C_n)\rightarrow \{1, 2, \dots, n\}$ be the vertex labeling function  given by
	
	$f(v_i) = \begin{cases}
	2i-1 &\mbox{if } 1\le i \le \lceil{\frac{n}{2}}\rceil\\
	2i-(n+1) &\mbox{if } \lceil{\frac{n}{2}}\rceil+1 \le i \le n~and~n~is~ odd\\
	2i-n &\mbox{if } \lceil{\frac{n}{2}}\rceil+1 \le i \le n~and~n~is~even 
	\end{cases}$
	
	This labeling first lists all the odd integers and then all the even integers. Hence, going from odd to even, exactly one edge receives the `$-$' sign and going from even to odd, exactly one more edge, i.e., $v_nv_1$ also receives `$-$' sign. For all other edges, we have the `+' sign. Hence, $\sigma^-(C_n)$=2. Thus the proof.
\end{proof}

Harary \cite{H53} defined the balanced signed graph as one in which every cycle has an even number of negative edges. From the above theorem, we have $\sigma^-(C_n)$=2. We find a link between parity signed cycle and balanced cycle. Every parity signed cycle has at least two negative edges. In fact, we have the following interesting result.

\begin{theorem}
	Every parity signed cycle C$_n$ is a balanced cycle.
\end{theorem} 
\begin{proof}
	\textbf{Theorem} \ref{negc} implies that negative homogeneous parity signed cycle is of even length. \textbf{Proposition} \ref{cyc2} gives us that for any cycle C$_n$ with $n$ vertices, $\sigma^-(C_n)$=2. Hence, $|E^-(C_n)|\geq 2$. 
	
	If $n$ is odd, then assigning the vertices from 1 to $n$, in the increasing order in a greedy assignment, we get every edge with vertices labeled with opposite parity, except a single edge with its vertices labeled 1 and $n$. Hence, in such a case, $|E^-(C_n)|=n-1$, which is even.
	
	In other cases, let us assume that the labeling begins from 1. When every time, the labels of the vertices switch from an odd integer to an even integer, there will be another switch from an even integer to an odd integer as explained above.
	
	Hence, we conclude that $|E^-(C_n)|$ is alway even. Therefore, 
	whether $n$ is odd or even, every parity signed cycle C$_n$ is a balanced cycle.
\end{proof}

The next result is about the importance of connectedness in parity signed graphs.

\begin{theorem}
	A positive homogeneous parity signed graph is not connected.
\end{theorem} 
\begin{proof}
	Assume that a positive homogeneous signed graph is a parity signed graph. Hence, labels of every pair of adjacent vertices are of the same parity. However, \textbf{Theorem} \ref{nont} guarantees that every connected non-trivial parity signed graph will have at least one negative edge. This forces us to think of a signed graph with at least two components where all vertices of one of the components are labeled with odd integers and all vertices of the other component are labeled with even integers. Therefore, the graph is convincingly not connected.
\end{proof}

\begin{proposition}\label{thm_star}
	For a star K$_{1,n}$ with $n+1$ vertices, $\sigma^-(K_{1,n})$=$\lceil{\frac{n}{2}}\rceil$.
\end{proposition}
\begin{proof}
	Let the vertices of the star K$_{1,n}$ be $\{v_i| i=1, 2, \dots, n+1\}$ and the edges be $\{v_1v_{i}| 2\le i \le n+1\}$. When $n$= 1 or 2, $\sigma^-(K_{1,n})$=1. Hence, let $n \ge 3$. 
	
	Assume that $n$ is odd. Then there are $n+1$ vertices and $n+1$ is even. Hence there are equal number of even and odd integers for labeling the vertices. If the non-pendant vertex, $v_1$, is given an odd label, then there will be ($\frac{n+1}{2}-1$) odd pendant vertices and $\frac{n+1}{2}$ even pendant vertices. Then the induced parity signed labeling implies that $|E^-(K_{1,n})|=\frac{n+1}{2}$. On the other hand, if the non-pendant vertex, $v_1$, is given an even label, then there will be $\frac{n+1}{2}-1$ even pendant vertices and $\frac{n+1}{2}$ odd pendant vertices. This also will give us $|E^-(K_{1,n})|=\frac{n+1}{2}$.
	
	Hence, if $n$ is odd, $\sigma^-(K_{1,n})=\frac{n+1}{2}$
	
	Now assume that $n$ is even. Then $n+1$ is odd and there are ($\frac{n}{2}+1$) odd integers and $\frac{n}{2}$ even integers. If the non-pendant vertex, $v_1$, is given an odd label, then there will be $\frac{n}{2}$ odd pendant vertices and $\frac{n}{2}$ even pendant vertices. Then the induced parity signed labeling implies that $|E^-(K_{1,n})|=\frac{n}{2}$. However, if the non-pendant vertex, $v_1$, is given an even label, then there will be ($\frac{n}{2}-1$) even pendant vertices and ($\frac{n}{2}+1$) odd pendant vertices. This will give us $|E^-(K_{1,n})|=\frac{n}{2}+1$.
	
	Hence, if $n$ is even, $\sigma^-(K_{1,n})=\frac{n}{2}$.
	
	From both the cases, we conclude that for a star K$_{1,n}$ with $n+1$ vertices, $\sigma^-(K_{1,n})$=$\lceil{\frac{n}{2}}\rceil$.
\end{proof}

The following easier proof is given to the Proposition \ref{thm_star} by the referee(s). Assume that the central vertex is odd. Then, for $n$ even, we get $-$ sign only for the edges with ends labeled with odd numbers. Their cardinality is $n/2$. For $n$ odd, we get $-$ sign only for the edges with ends labeled with even numbers. Their cardinality is $n/2$. If the label of the central vertex is even, then the argument is similar.
\begin{proposition}\label{compl}
	For a complete graph K$_n$ with $n \ge 2$ vertices, $\sigma^-(K_n) = \lfloor{\frac{n}{2}}\rfloor \lceil{\frac{n}{2}}\rceil$.
\end{proposition} 
\begin{proof}
	Let K$_n$ be a complete graph  with $n$ vertices. In any complete graph, every pair of vertices is adjacent. Let the vertices be $\{v_i| i=1, 2, \dots, n\}$. Let $f:V(K_n)\rightarrow \{1, 2, \dots, n\}$ be the vertex labeling function given by $f(v_i)=i, i=1, 2, \dots, n$. This obviously, gives a parity signed labeling of the edges of K$_n$.
	%nparagraph break
	Let A and B be set of all vertices of K$_n$ labeled with odd and even postive integers, respectively. Hence, every edge of K$_n$ between vertices of A will get `+' in the induced parity signed labeling. Similarly, every edge of K$_n$ between vertices of B will also get `+' in the induced parity signed labeling.
	
	Consequently, all edges of K$_n$ between vertices of A and B will get `$-$' labels. The number of edges of K$_n$ between vertices of A and B is $|A|.|B|$.
	
	If $n$ is even, then $|A|=|B|=\frac{n}{2}$. Hence, $\sigma^-(K_n) = \frac{n^2}{4}$.
	
	If $n$ is odd, then $|A|=\frac{n+1}{2}$ and $|B|=\frac{n-1}{2}$.
	
	Hence, $\sigma^-(K_n) = (\frac{n+1}{2})(\frac{n-1}{2})$=$\lceil{\frac{n}{2}}\rceil \lfloor{\frac{n}{2}}\rfloor$.
	
	Thus we conclude that for K$_n$, $n \ge 2$, $\sigma^-(K_n)$=$\lfloor{\frac{n}{2}}\rfloor \lceil{\frac{n}{2}}\rceil$.
\end{proof}
Please note that $\sigma^-(K_n)$ is nothing but the product of the number of odd integers and the number of even integers of $\{1, 2, 3, \dots \}$. We thank the referee(s) for this elegant observation.
\section{Some Characterizations}

\begin{theorem}\label{Th_St}
	Let $T$ be a parity signed tree. Then $\sigma^-(T)=|E^-(T)|$ if and only if $T$ is $K_{1,n}$ where $n\in \mathbb{N}$ is odd.
\end{theorem}
\begin{proof}
	Assume that the parity signed tree is $K_{1,n}$ where $n$ is odd. Then, by Proposition \ref{thm_star}, $\sigma^-(K_{1,n})=\frac{(n+1)}{2}$. $K_{1,n}$ has $n+1$ vertices which is an even number. If the central vertex has an odd integer as label, then there are exactly $\frac{(n+1)}{2}$ pendant vertices with even integer labels, by which we get $|E^-(T)|=\frac{(n+1)}{2}$. If the central vertex is labeled with an even integer, then there are exactly $\frac{(n+1)}{2}$ pendant vertices with odd integer labels, by which we get $|E^-(T)|=\frac{(n+1)}{2}$. Hence, $\sigma^-(K_{1,n})=|E^-(K_{1,n})|$, when $n$ is odd.
	
	On the contrary, assume that the parity signed tree is not $K_{1,n}$ where $n\in \mathbb{N}$ is odd.\\
	\textbf{Case 1}:  The parity signed tree is $K_{1,n}$ where $n\in \mathbb{N}$ is even.
	
	Here, the total number of vertices $n+1$ is odd. Therefore, there are $\frac{n}{2}$ even integers and $\frac{n}{2}+1$ odd integers from 1 to $n+1$. Hence, if the central vertex is odd, then $|E^-(K_{1,n})|=\frac{n}{2}$ and if the central vertex is even, then $|E^-(K_{1,n})|=\frac{n}{2}+1$. i.e., the \textit{rna} is not equal to the negative edges of $T$.\\
	\textbf{Case 2}:  The parity signed tree is not a star.
	
	Here, the tree must have at least one path $P_4$ of length having a pendant vertex. Let the vertices of this $P_4$ be  $v_1, v_2, v_3$ and $v_4$ with $v_1$ being the pendant vertex and other being the consecutive vertices. We now show that by a swapping of the labels of two adjacent vertices, the number of negative edges vary. If the tree itself is $P_4$, the result is obvious. Hence, we look for non-obvious cases.\\
	\textbf{Subcase 2. 1}: Assume that $v_2$ is of degree two. We can assign the labels 1, 2, 3 and 4 to $v_1, v_2, v_3$ and $v_4$, respectively. This provides 4 negative edges to the signed graph. Now, we swap the labels 1 and 2. Now, labels 1 and 3 are adjacent. Hence, this new labeling reduces the negative edges of the signed graph by 1.\\
	\textbf{Subcase 2. 2}: Assume that the vertex $v_2$ is of even degree greater than 2. i.e., the vertex $v_2$ has some $2k, k\in \mathbb{N}$ neighbours other than $v_1$ and $v_3$. After assigning the labels 1, 2, 3 and 4 to $v_1, v_2, v_3$ and $v_4$, respectively, we assign labels 5, 6, $\dots$ in that order to the $2k$ vertices adjacent to the vertex $v_2$. Hence, the parity signed graph has at least 3+$k$ negative edges. Now, swap the labels 1 and 2. Now, labels 1 and 3 are adjacent. Hence, this new labeling reduces the negative edges of the signed graph by 1.\\
	\textbf{Subcase 2. 3}: Assume that the vertex $v_2$ is of odd degree
	greater than 2. i.e., the vertex $v_2$ has some $2k+1, k\in \mathbb{N}\cup\{0\}$ neighbours other than $v_1$ and $v_3$. After assigning the labels 1, 2, 3 and 4 to $v_1, v_2, v_3$ and $v_4$, respectively, we assign labels 5, $\dots$ in that order to the $2k+1$ vertices adjacent to the vertex $v_2$. Hence, the parity signed graph has at least 3+$k$+1 negative edges. Now, swap the labels 1 and 2. Now, labels 1 and 3 are adjacent. Hence, in the $P_4$, two edges are negative and one edge is positive. Further, the number of negative edges from $v_2$ to the neighbours other than $v_1$ and $v_2$ is also reduced by 1.
	
	These three subcases only are needed as we deal with only the edges incident to the vertex $v_2$. In fact, what we have shown is that number of negative edges of the parity signed tree is dependent on the labels of its vertices.
	
	To sum up from these two cases, we conclude that depending on the labeling of the vertices of the tree, the negative edges vary if the parity signed tree is not $K_{1,n}$ where $n\in \mathbb{N}$ is odd. Thus the result.
\end{proof}

%%Corona Product of Graphs%%
\begin{defn}\cite{W}
	Let $J$ and $K$ be any two graphs with $j$ and $k$ vertices, respectively. The corona graph, $J\circ K$, is a graph formed by one copy of $J$ and $j$ copies of $K$ by making every vertex of a $K$ adjacent to exactly one vertex of $J$.
\end{defn}

\begin{theorem}\label{negh}
	The corona, $C_n\circ K_1$, is negative homogeneous parity signed graph if and only if $n$ is even.
\end{theorem} 
\begin{proof}
	Let the corona, $C_n\circ K_1$, be negative homogeneous. This means that every edge of $C_n$ and every pendant edge is negative.
	
	To have every edge on the cycle negative, every pair of adjacent vertices must have labels of opposite parity. However, $C_n\circ K_1$ has 2$n$ vertices. Therefore, only if $n$ is even, every pair of adjacent vertices on the cycle $C_n$ will have labels of opposite parity and labels of vertices of pendent edges will also have labels of opposite parity. 
	
	Conversely, assume that $n$ is even.
	Let the vertices of the cycle be $u_1, u_2, \dots, u_n$ where $u_1u_2, u_2u_3,$ $\dots,$ $u_{n-1}u_n$ and $u_nu_1$ are its edges. Let the pendant vertex corresponding to the vertex $u_i$ be $v_i$ for $1\leq i \leq n$. Then, the labeling $u_i \rightarrow i$ and $v_i \rightarrow (2n+1-i)$ give the desired labeling. 
\end{proof}
\begin{theorem}
	The corona, $K_n\circ K_1$, is negative homogeneous parity signed graph if and only if $n \le 2$.
\end{theorem} 
\begin{proof}
	When $n$=1 or 2, the result is obvious. When $n$=3, $K_n$ is nothing but $C_n$ and $K_3\circ K_1$ is not a negative homogeneous parity signed graph (\textbf{Theorem} \ref{negh}). Hence, let $n\geq4$.
	
	Now, $K_n\circ K_1$ has $K_n$ as a subsignedgraph.  As $n\geq4$, $K_n$ has $C_3$ as a subsignedgraph which cannot be negative homogeneous. Hence, $K_n\circ K_1$ is not a negative homogeneous parity signed graph if $n > 2$.
\end{proof}
\begin{theorem}
	Let $S$ be a signed graph formed by two negative homogeneous parity signed graphs connected by a bridge. Then $S$ is a parity signed graph.
\end{theorem}
\begin{proof}
	Given that the signed graph $S$ is formed by two negative homogeneous parity signed graphs connected by a bridge. Let the two negative homogeneous signed graphs and the bridge be $S_1, S_2$ and $uv$, respectively. Let $u\in V(S_1)$ and $v\in V(S_2)$. Since $S_1$ and $S_2$ are negative homogeneous, both of them will have vertices labeled with odd and even integers. Taking vertices $u$ and $v$ having labels of the same (opposite) parity and joining them by positive (negative) edge  produces the desired result.
\end{proof}

\section{Conclusion}

We have defined a parity signed graph, its `rna' and `adhika' numbers. We have evaluated the `rna' number of path, star, cycle and complete graph. The balanced structure of a parity signed cycle is also explored. A couple of characterizations associated with complete graphs and some other graphs are also found out. Characterization of parity signed graphs among many other classes of signed graphs is an open problem. In addition to this, there is ample scope for research in finding the relations between the `rna' number and the `adhika' number of signed graphs and counting the number of possible parity signed labelings of a signed graph. Proposition \ref{compl} and Theorem \ref{Th_St} motivate us to propose the following conjecture.
\begin{con}
	Let $G$ be a parity signed graph. Then $\sigma^-(G)=|E^-(G)|$ if and only if $G$ is either $K_{1,n}$, $n$ odd or $K_{n}$.
\end{con}

\section*{Acknowledgments}
We dedicate this paper to the inspirational personality of Dr. B. Devadas Acharya who advanced human knowledge in the areas of Signed Graphs, Domination Theory and Hypergraphs.

\end{document}